\documentclass[a4paper,11pt]{amsart}
\usepackage{extarrows}
\usepackage{amssymb}
\usepackage{mathrsfs}
\usepackage{amsmath,amssymb,amsthm,latexsym,amscd,mathrsfs}
\usepackage{indentfirst}
\usepackage{extarrows}
\usepackage[all]{xy}
\usepackage{graphicx}
\usepackage{multirow}

\numberwithin{equation}{section} \theoremstyle{plain}
\newtheorem{thm}{Theorem}[section]

\newtheorem{conj}[thm]{Conjecture}
\newtheorem{exmp}[thm]{Example}
\newtheorem{prob}[thm]{Problem}

\newtheorem{ques}[thm]{Question}

\newtheorem{ack}{Acknowledgements}   
 \makeatletter

\newcommand{\fol}{\mathcal{F}}
\newcommand{\sphere}{\mathbb{S}}
\DeclareMathOperator{\Diff}{Diff}
\DeclareMathOperator{\Isom}{Isom}

\makeatother
\title[Differential topology v.s isoparametric foliations]{Differential topology interacts with isoparametric foliations}
\author[J.Q. Ge]{Jianquan Ge}
\address{School of Mathematical Sciences, Laboratory of Mathematics and Complex Systems, Beijing Normal
University, Beijing 100875, P.R. CHINA.}
\email{jqge@bnu.edu.cn}
\author[C. Qian]{Chao Qian}\address{School of Mathematics and Statistics, Beijing Institute of Technology, Beijing 100081, P.R.China. }\email{qianchao\_1986@163.com}
\thanks{The second author is the corresponding author.}
\date{}
\thanks{The project is partially supported by the NSFC (No. 11331002 and No. 11401560) and the Fundamental Research
Funds for the Central Universities. }
\begin{document}
\maketitle

%%%%%%%%%%%%%%%%%%%%%
\begin{abstract}
In this note, we discuss the interactions between differential topology and isoparametric foliations, surveying some recent progress and open problems.
\end{abstract}
%%%%%%%%%%%%%%%%%%%%%%%

%%%%%%%%%%%%%%%%%%%%%%
\section{Introduction}\label{introduction}
As is well known, two of the main topics that differential topology studies are smooth structures on manifolds and smooth mappings between them. Since the surprising discovery of exotic spheres by Milnor \cite{Mil56} in 1956, existence and non-existence of exotic smooth structures have obtained worldwide attention and highly intensive study. Recall that an exotic sphere is a closed smooth manifold which is homeomorphic but not diffeomorphic to the unit sphere $\mathbb{S}^n$. In 1963, Kervaire and Milnor \cite{KM63} gave a detailed investigation of the group $\Theta_n$ of h-cobordism classes of oriented homotopy $n$-spheres. Here a homotopy $n$-sphere $\Sigma^n$ is a closed smooth manifold which has the homotopy type of $\mathbb{S}^n$. In this case, $\Sigma^n$ is known to be homeomorphic to $\mathbb{S}^n$ (cf. \cite{JW08}).
It is well known that $\Theta_n$ is isomorphic to
$\Gamma_n$, the group of oriented twisted $n$-spheres.
What is more, according to Cerf \cite{Ce70}, $\Theta_n$ is isomorphic to the mapping class group $\pi_0\mathrm{Diff}^+(\mathbb{S}^{n-1})$ by $\pi_0\mathrm{Diff}^+(\mathbb{S}^{n-1})\rightarrow \Theta_n=\Gamma_n, [\phi]\mapsto \Sigma_{\phi}:=D^n\cup_{\phi}D^n$. Note that
$\Sigma_{\phi}$ depends only on the isotopy class of $\phi \in \mathrm{Diff}^+(\mathbb{S}^{n-1})$. Motivated by $\Theta_n$, two inertia
groups $I_0(M)$ and $I_1(M)$ are defined for a closed oriented manifold $M$ (cf. \cite{Le70}). In fact,
$I_0(M)$ is related to the study of exotic smooth structures on $M$, and
$I_1(M)$ contributes to the group $\Gamma(M)$ of pseudo-isotopy classes of diffeomorphisms of $M$.

On the one hand, from the local viewpoint of Riemannian geometry, one of the central problems is to determine the classes of manifolds with special curvature properties, for instance, manifolds with positive/nonnegative sectional, Ricci or scalar curvature. Therefore, the curvature properties of manifolds with exotic smooth structures, especially exotic spheres, are very interesting.
%As pointed by Grove and Ziller (cf. \cite{GZ00}), since Milnor's
%discovery of exotic spheres, one of the most intriguing problems in Riemannian %Geometry has been whether there are exotic spheres with positive sectional %curvature.
We refer to \cite{GZ00} and \cite{JW08} for more details and the progress of this subject.

On the other hand, from the global viewpoint, it is rather fascinating to study singular Riemannian foliations on manifolds, especially isoparametric foliations, which are geometric generalizations of manifolds with isometry group actions. It was E. Cartan who firstly gave a thorough study of isoparametric foliations on the unit spheres in 1930's. Up to now, the study
of isoparametric foliations has become a highly influential field in differential geometry. We recommend \cite{Ce08}, \cite{Ch12}, \cite{GT12}, \cite{QT14}, \cite{TXY14}, \cite{TY13}, \cite{TY15} and  \cite{Th00} for a systematic and complete survey of isoparametric foliations and their applications.

The contents of this note are mainly extracted from \cite{Ge14} and \cite{GR13}, and are organized as follows. In Section 2 we will recall the basic notations, examples and the fundamental facts that lead to the interaction between differential topology and isoparametric foliations. In Section 3 we will discuss the interaction in dimension four which in particular produces a classification of $4$-manifolds with singular Riemannian foliations. A related conjecture of K. Grove about nonnegative curvature is also discussed. In Section 4 we will pay attention to isoparametric foliations on homotopy spheres and some relations with diffeomorphism groups, for example $\pi_0\mathrm{Diff}^+(\mathbb{S}^{4})$. In Section 5 some applications to exotic smooth structures and inertial groups will be presented.

\section{Singular Riemannian foliations and isoparametric foliations}
In this section, we will recall the basic notations and results on singular Riemannian foliations (cf. \cite{Mol88, Wa87, GW09, ABT13, Th10}).
Let $M$ be a complete Riemannian manifold. A \emph{transnormal system} $\fol$ is a decomposition of $M$ into complete connected injectively immersed submanifolds, called \emph{leaves}, such that each geodesic emanating perpendicularly to one leaf remains perpendicular to the leaves at all its points.
A \emph{singular Riemannian foliation} (SRF) is a transnormal system $\fol$ which
is also a \emph{singular foliation}, i.e., such that there are smooth vector
fields $X_i$ on $M$ that span the tangent spaces $T_pL_p$ to the leaf $L_p$ through each point $p\in M$. A SRF $\fol$ is called a \emph{polar foliation} if through each point of $M$ passes a complete immersed submanifold, called a section, that intersects all the leaves of $\fol$ orthogonally. A leaf of maximal dimension is called a \emph{regular} leaf, and its codimension is defined to be the \emph{codimension} of $\fol$. Leaves of lower dimension are called singular leaves. In particular, a singular Riemannian foliation all of whose leaves have the same dimension is called a (regular) Riemannian foliation.

A singular Riemannian foliation $(M, \fol)$ of codimension 1 is called an \emph{isoparametric} foliation if all the regular leaves have constant mean curvature. Each regular leaf of an isoparametric
foliation is called an \emph{isoparametric hypersurface}, and the singular leaves are called \emph{focal submanifolds}. Moreover, if each regular leaf has constant principal curvatures, then $\fol$ is called a \emph{totally} isoparametric foliation (cf. \cite{GTY14}). According to \cite{Wa87}, there is essentially a correspondence between transnormal (resp. isoparametric) functions on $M$ and transnormal
systems (resp. isoparametric foliations) of codimension 1 on $M$.

Next to the basic notations, we are in a position to provide some interesting examples of SRF.
\begin{exmp}
\item[(1).] Homogeneous SRF:

Let $G$ be a Lie group that acts on a Riemannian manifold $M$ by isometries, and let
$\fol$ be the partition of $M$ by the orbits of $G$. Then $\fol$ is a singular Riemannian foliation on $M$, called homogeneous SRF. If the codimension of $\fol$ is equal to 1, i.e., the isometric action is a cohomogeneity one action, then $\fol$ is a homogeneous isoparametric foliation. Homogeneous isoparametric foliations are totally isoparametric.

\item[(2).] Isoparametric hypersurfaces in unit spheres (\cite{Th00, Ce08, Ch12}):

Due to E. Cartan, isoparametric hypersurfaces in real space forms have constant principal curvatures, and the spherical case is more complicated. Let $g$ be the number of distinct principal curvatures. For isoparametric hypersurfaces in $\mathbb{S}^n$, $g$ must be $1,2,3,4$ or $6$. If g=$1,2,3$ or $6$,
isoparametric hypersurfaces are all homogeneous (cf. \cite{Mi13}). For the case $g=4$, isoparametric hypersurfaces must be OT-FKM type or homogeneous except for one case (cf. \cite{CCJ07, Ch13}).
Consequently, there exist infinitely many inhomogeneous, but totally isoparametric foliations on unit spheres.

\item[(3).] Isoparametric foliations in compact symmetric spaces:

In the complex projective space $\mathbb{C}P^n$ with the Fubini-Study
metric, every totally isoparametric foliation must be homogeneous (cf. \cite{GTY14}). Moreover, in $\mathbb{C}P^{2m}$ all isoparametric foliations are homogeneous and in each $\mathbb{C}P^{2m+1}$ there exist inhomogeneous isoparametric foliations.
Based on the classification result for $\mathbb{S}^n$ and a classification of complex structures, \cite{DV12} obtained the classification of isoparametric foliations on all $\mathbb{C}P^n$ except for one case. Similarly, \cite{DVG15} obtained the classification of isoparametric foliations on all $\mathbb{H}P^n$ except for one case.  So far,
the classification problem of isoparametric foliations on compact symmetric spaces is far from being touched, except for the cases $\mathbb{S}^n$, $\mathbb{C}P^n$ and $\mathbb{H}P^n$ .

\item[(4).] SRF and Regular Riemannian foliation on (homotopy) spheres:

It was shown by Thorbergsson \cite{Th91} that polar foliations of high codimension ($\geq2$) on unit spheres are homogeneous. However, non-polar SRF (of high codimension) can be inhomogeneous and are abundant on unit spheres (cf. \cite{Ra14}).
Compared to SRF, regular Riemannian foliations are rather rare on unit spheres, even on homotopy spheres. They occur only when the dimension of the leaves is $1$, $3$ or $7$ (cf. \cite{LW13}). Nevertheless, they have not yet been classified on homotopy spheres up to foliated diffeomorphisms.

\end{exmp}
Motivated by the examples (2-3) above, it is natural to pose the following
\begin{ques}
Is every totally isoparametric foliation on a compact symmetric space other than $\mathbb{S}^n$
a homogeneous isoparametric foliation?
\end{ques}

Now, we turn to the topological side of SRF with codimension 1. In this note, by a \emph{foliated diffeomorphism} between two foliated manifolds, we mean a diffeomorphism maps leaves to leaves, and such foliations are called \emph{equivalent}. And two equivalent foliated manifolds are denoted by $(M,\fol)\cong(M',\fol')$. Note that the foliated diffeomorphism here needs not to be an isometry.

According to \cite{Mol88}, a codimension 1 SRF $\fol$ on a closed simply connected manifold $N$ has exactly two closed singular leaves $M_{\pm}$, and $N$ has a double disk bundle decomposition (DDBD) by two disk bundles $E_{\pm}$ of the normal vector bundles $\xi_{\pm}$ over $M_{\pm}$ of rank $m_{\pm}>1$, i.e., $N\cong E_{\varphi} := E_+\cup_{\varphi}E_-$ with the gluing diffeomorphism $\varphi: \partial E_+\rightarrow \partial E_-$ (cf. \cite{Ge14}). The codimension 1 SRF $\fol$ on $N$ induces a codimension 1 SRF $\fol_{\varphi}$ on $E_\varphi$. The leaves of $\fol_{\varphi}$ are just the concentric tubes around the zero sections in $E_{\pm}$. Hence, the equivalence class of $(N, \fol)$ can be represented by $(E_{\varphi}, \fol_{\varphi})$.

Conversely, given a DDBD, i.e., given disk bundles $E_{\pm}$ over $M_{\pm}$ with rank greater than 1 and a diffeomorphism  $\varphi: \partial E_+\rightarrow \partial E_-$, we let $\fol_{\varphi}$ be the singular foliation consisting of concentric tubes on $E_{\varphi}:=E_+\cup_{\varphi} E_-$. By the fundamental construction theorem in \cite{QT13}, there exists a bundle-like metric $g_{\varphi}$ such that $(E_{\varphi}, \fol_{\varphi})$ becomes an isoparametric foliation.

It follows that isoparametric foliations on closed simply connected manifolds require no more on the topology than codimension 1 SRF, and they are all equivalent to the topological condition DDBD structure. Moreover, to study the classification of equivalence classes of isoparametric foliations, one only needs to study foliations in the form $(E_{\varphi}, \fol_{\varphi})$ determined by two pairs of disk bundles $E_{\pm}\subset \xi_{\pm}$ and gluing diffeomorphisms $\varphi:\partial E_+\rightarrow \partial E_-$. In this way, differential topology comes up to take place of submanifold geometry in the theory of isoparametric foliations, once we are only concerned and satisfied with a classification up to foliated diffeomorphisms other than ambient isometries. This way has its priority over the classical way, in the sense that now the existence problem disappears if provided with a DDBD structure (since now foliations are not sensitive about metric), while it is extremely difficult to find an example under a general fixed metric. In particular, in low dimensional cases the topological condition DDBD structure suffices to determine the differentiable structure as illustrated in the following section.

\section{Interaction between $4$-manifolds and SRF}
In this section, we will be concerned with the singular Riemannian foliations on simply connected $4$-manifolds. Geometry and topology of $4$-manifolds form an extremely rich and also complicated research field. This is the lowest dimension in which exotic smooth structures arise, e.g., the noncompact $4$-spaces $\mathbb{R}^4$ and compact $m\mathbb{CP}^2$\#$ n\overline{\mathbb{CP}^2}$ for many pairs of ($m\geq 1, n\geq 2$) (cf. \cite{FQ90, AP2}). Up to now, it is not known
whether there is a $4$-manifold with only one smooth structure, even for $\mathbb{S}^4$, the affirmative side of which is called the \emph{smooth Poincar\'{e} conjecture}. Therefore, it is natural to study $4$-manifolds with additional structures, especially $4$-manifolds with SRF. Some partial classification results
are listed as the following.
\begin{thm} Known classification results of closed $1$-connected $4$-manifolds with:
\begin{itemize}
\item[(1)] Homogeneous SRF:
\item[a)] If $N$ admits a cohomogeneity one action, then $N$ is diffeomorphic to one of the $4$ standard manifolds $\mathbb{S}^4$, $\mathbb{CP}^2$, $\mathbb{S}^2\times\mathbb{S}^2$, or $\mathbb{CP}^2\#- \mathbb{CP}^2$ (\cite{Pa86, GH87, GZ00, Ho10})
\item[b)] If $N$ admits a $T^2$ action (with cohomogeneity two),  then $N$ is diffeomorphic to a connected sum of copies of standard $\sphere^4$, $\pm\mathbb{CP}^2$ and $\sphere^2\times\sphere^2$ (\cite{OR70, OR74}).
\item[b)] If $N$ admits a $\sphere^1$ action (with cohomogeneity three), then $N$ is diffeomorphic to a connected sum of copies of standard $\sphere^4$, $\pm\mathbb{CP}^2$ and $\sphere^2\times\sphere^2$ (\cite{Fi77, Fi78}).
\item[(2)] General SRF:
\item[a)] If $N$ is a homotopy $4$-sphere admitting a codimension 1 SRF, then $N\cong\sphere^4$ (\cite{GT13}).
\item[b)] If $N$ admits a SRF $\fol$ of codimension 3, then $\fol$ must be a homogeneous SFR induced by some $\sphere^1$ action (\cite{GaRa13}).
\end{itemize}
\end{thm}
By analyzing the topological restrictions (e.g., cohomology, sphere bundles over spheres, isotopy classes of diffeomorphisms, etc.) given by the DDBD structure, \cite{GR13} was able to give a complete differentiable classification of $4$-manifolds with isoparametric foliations and also a classification of isoparametric foliations up to equivalence classes in dimension $4$. Moreover, for SRF of higher codimension \cite{GR13} also obtained the differentiable classification by some inductive method with respect to the number of edges of the leaf space. These generalize all the results mentioned in the theorem above. Explicitly, \cite{GR13} showed
\begin{thm}\label{thm4mfds}(\cite{GR13}) Let $N$ be a closed simply connected $4$-manifold.
\begin{itemize}
\item[(a)] If $N$ admits a SRF of codimension $1$, i.e., a DDBD structure, then $N$ is diffeomorphic to one of the $5$ standard manifolds $\mathbb{S}^4$, $\mathbb{CP}^2$, $\mathbb{S}^2\times\mathbb{S}^2$, or $\mathbb{CP}^2\#\pm \mathbb{CP}^2$.
\item[(b)] If $N$ admits a SRF of codimension greater than $1$, then $N$ is diffeomorphic to a connected sum of copies of standard $\sphere^4$, $\pm\mathbb{CP}^2$ and $\sphere^2\times\sphere^2$.
\end{itemize}
\end{thm}

\begin{thm}(\cite{GR13})
Let $N$ be a closed simply connected $4$-manifold admitting a SRF $\fol$ of codimension $1$ with regular leaf $M$ and two singular leaves $M_{\pm}$. Then the equivalence class of $(N, \fol)$ can be uniquely represented in terms of $(M, M_{\pm})$
and classified in the following Table $1$:
\begin{table}[!htb]%[!hbp]
\caption{SRF of codim 1 on closed simply connected $4$-manifolds}\label{table-SRF of codim1}
\begin{tabular}{|c|c|c|c|c|c|}
\hline
\multirow{2}{*}{$N$} & \multicolumn{2}{|c|}{$\fol$} & \multicolumn{3}{|c|}{Properties} \\
\cline{2-6}
& $M$ & $M_{\pm}$ & Homog & T-Isopar & Isopar\\
\hline
\multirow{3}{*}{$\sphere^4$}  & $L(1,1)$ &  $pt$, $pt$ &\multirow{6}{*}{Yes} & \multirow{9}{*}{Yes} & \multirow{10}{*}{Yes}\\
\cline{2-3}
& $L(0,1)$ &  $\sphere^1$, $\sphere^2$ &&&\\
\cline{2-3}
& $SO(3)/(\mathbb{Z}_2\oplus\mathbb{Z}_2)$ &  $\mathbb{RP}^2$, $\mathbb{RP}^2$ &&&\\
\cline{2-3}
\cline{1-3}
\multirow{2}{*}{$\mathbb{CP}^2$}  & $L(1,1)$ &  $pt$, $\sphere^2$ & & &\\
\cline{2-3}
& $L(4,1)$ &  $\mathbb{RP}^2$, $\sphere^2$ &&&\\
\cline{2-3}
\cline{1-3}
$\sphere^2\times\sphere^2$ & $L(2m,1)$,   $m\geq 0$ &  $\sphere^2$, $\sphere^2$ &  & &\\
\cline{1-4}
\multirow{2}{*}{$\mathbb{CP}^2\# \mathbb{CP}^2$}  & $L(1,1)$ &  $\sphere^2$, $\sphere^2$ &\multirow{2}{*}{No} &&\\
\cline{2-3}
& $L(2,1)$ &  $\sphere^2$, $\sphere^2$ &&&\\
\cline{2-3}
\cline{1-4}
\multirow{2}{*}{$\mathbb{CP}^2\# -\mathbb{CP}^2$}  & $L(2m+1,1)$,   $m\geq0$  &  $\sphere^2$, $\sphere^2$ & Yes &&\\
\cline{2-5}
& $L(0,1)$ &  $\sphere^2$, $\sphere^2$ & No & Unknown &\\
\cline{2-3}
\hline
\end{tabular}
\end{table}\\
where the column ``Homog" (resp. ``T-Isopar", ``Isopar") means whether there exist a homogeneous (resp. totally isoparametric, isoparametric) representative in the foliated diffeomorphism class.

\end{thm}

To conclude this section, we state a conjecture of K. Grove which relates nonnegative sectional curvature with isoparametric foliation.
\begin{conj}(\cite{GR13})
Any closed simply connected non-negatively curved manifold admits a DDBD structure, i.e., a SRF of codimension $1$ (under some metric).
\end{conj}

By Theorem \ref{thm4mfds} (a),  an affirmative answer to the conjecture above will solve
\begin{conj}\label{Bott}(\cite{GR13,GH87})
A closed simply connected non-negatively curved $4$-manifold is diffeomorphic to one of the $5$ standard manifolds $\mathbb{S}^4$, $\mathbb{CP}^2$, $\mathbb{S}^2\times\mathbb{S}^2$, or $\mathbb{CP}^2\#\pm \mathbb{CP}^2$.
\end{conj}
Note that Conjecture \ref{Bott} is still open even for the homeomorphism case, whereas it is indeed true even up to equivariant diffeomorphism if
provided with an isometric circle action (see e.g. \cite{GW14} and references therein).

\section{Isoparametric foliations on homotopy spheres}
The study of codimension $1$ SRF and isoparametric foliations on homotopy spheres was
initiated by Ge and Tang in \cite{GT13}. Based on their research, they proposed the following problem.
\begin{prob}\label{prob}(\cite{GT13})
Are there always isoparametric foliations on $\Sigma^n$ ($n\neq4$) with the same focal submanifolds as those on $\sphere^n$?
\end{prob}
To attack this problem, Qian and Tang in \cite{QT13} obtained a fundamental construction: a bundle-like metric such that a DDBD structure is isoparametric under this metric, as mentioned in Section 2. This established a bridge between differential topology and the theory of isoparametric foliations. As an immediate corollary, each homotopy sphere $\Sigma^n$ ($n>4$) admits an isoparametric foliation with two points as the focal submanifolds since in dimension greater than $4$, a homotopy sphere is always a twisted sphere (a celebrated result of Smale), i.e., $\Sigma_{\phi}^n=D^n\cup_{\phi}D^n$ for some $\phi\in \Diff^+(\sphere^{n-1})$. This answer Problem \ref{prob} partially.

Following this way, to answer Problem \ref{prob} completely it suffices to find the DDBD structure for $\Sigma^n$ with the same disk bundles as those for $\sphere^n$. This was completed in \cite{Ge14}, where even more were obtained. In fact, given a DDBD structure on $\sphere^n$, say $\sphere^n\cong E_{\varphi}=E_+\cup_{\varphi}E_-$, there is a DDBD sturcture on $\Sigma_{\phi}^n$ in the form $\Sigma_{\phi}^n\cong E_{d_{\phi}\circ\varphi}=E_+\cup_{d_{\phi}\circ\varphi}E_-$, and vice versa. Here $d_{\phi}\in\Diff(\partial E_-)$ is determined by $\phi\in\Diff^+(\sphere^{n-1})$. Moreover, this correspondence preserves equivalence classes. In conclusion, we have
\begin{thm}(\cite{Ge14})
 For each homotopy sphere $\Sigma^n$ ($n>4$), there exists a 1-1 correspondence between the sets of equivalence classes of isoparametric foliations on $\Sigma^n$ and $\sphere^n$. Moreover, the disk bundles of the corresponding DDBDs coincide with each other.
\end{thm}
Coarsely speaking, this theorem tells us that each homotopy sphere $\Sigma^n$ ($n>4$) is not only a twisted sphere but also a ``DDBD-twisted" sphere.

Consequently, the classification of equivalence classes of isoparametric foliations on homotopy spheres is equivalent to the case on $\sphere^n$. However, this classification on $\sphere^n$ ($n>4$) is still far from reached though almost completed for the round metric.

Considering the classification problem in $N=\mathbb{S}^n$ for example, in the following we introduce the three steps proposed in \cite{Ge14}.

\textbf{\underline{Step 1.}} Classify disk bundles $E_{\pm}$, s.t., $N=E_{\varphi}=E_+\bigcup_{\varphi}E_-$.

\textbf{\underline{Step 2.}} Classify isotopy classes of such $\varphi$, i.e., compute the subset
$$G_{N}(E_{\pm}):=\Big\{[\varphi]\in\pi_0\Big(\Diff(\partial E_+\rightarrow \partial E_-)\Big)\mid N\cong E_{\varphi}\Big\}$$

\textbf{\underline{Step 3.}} Compute the action:
$$\begin{array}{cccc}
\beta:& \pi_0(\Isom_b(\partial E_{\pm}))\times G_{N}(E_{\pm})& \rightarrow& G_{N}(E_{\pm})  \\
&([f_{\pm}],\quad [\varphi])&\mapsto& [f_{-}\circ\varphi\circ f_{+}^{-1}]
\end{array}$$
For each pair $E_{\pm}$ of Step 1, we have exactly $|G_{N}(E_{\pm})/\beta|$ equivalence classes of isoparametric foliations on $N$.

Before going on we introduce some notations. Let $\Isom_b(E)$ be smooth self bundle maps of $E$ preserving Euclidean metric, $\Isom_b(\partial E)$ be the restriction of $\Isom_b(E)$ to $\partial E$, $\Diff_E(\partial E)\subset\Diff(\partial E)$ be the subgroup of diffeomorphisms of $\partial E$ that extendable to $E$. The following result ensures the feasibility of Step 3.
\begin{thm}(\cite{Ge14})
For $\varphi_i:\partial E_{+}\rightarrow\partial E_{-}$, $i=0,1$,  $$(E_{\varphi_0},\fol_{\varphi_0})\cong(E_{\varphi_1},\fol_{\varphi_1})\Longleftrightarrow \exists f_{\pm}\in \Isom_b (\partial E_{\pm})~~s.t.~~[\varphi_1]=[f_{-}\circ\varphi_0\circ f_{+}^{-1}].$$
\end{thm}

In general, the set $G_{N}(E_{\pm})$ in Step 2 is intriguing, since we do not know whether there exists a diffeomorphism $\Phi: E_{\varphi_0}\rightarrow E_{\varphi_1}$ such that $\Phi(E_{\pm})=E_{\pm}$.
However, when one of $E_{\pm}$ is a disk $D^n$, this is indeed true, thanks to the Disk Theorem of Palais.
In this case, it turns out that $G_{N}(E_{\pm})=\Upsilon_{\varphi}$, a subset of $G_{N}(E_{\pm})$ defined as follows which is relatively easier to treat with.
$$\Upsilon_{\varphi}:=\Big\{[h_{-}\circ\varphi\circ h_{+}^{-1}]\mid [h_{\pm}]\in \pi_0(\Diff_{E_{\pm}}(\partial E_{\pm}))\Big\}\subset G_{N}(E_{\pm}), \quad \forall~~\varphi\in G_{N}(E_{\pm}).$$
Obviously $\Upsilon_{\varphi}$ is just an orbit of $\beta$ extended to $\pi_0(\Diff_{E_{\pm}}(\partial E_{\pm}))$ on $G_{N}(E_{\pm})$.
There is yet a smaller subset $\Gamma_{\varphi}:=\Big\{[\varphi']\mid \varphi'~\textit{is pseudo-isotopic to}~\varphi\Big\}\subset\Upsilon_{\varphi}$.

Combining these with known facts about $\pi_0(\Diff_{D^n}(\sphere^{n-1}))$, it follows that$$G_{\mathbb{S}^n}^+(D^n,D^n)=\Upsilon_{\varphi}^+=\Gamma_{\varphi}^+=\{[\varphi]\}, \quad n\neq 5.~~\Rightarrow |G_{\mathbb{S}^n}(D^n,D^n)/\beta|=1;$$
$$\pi_0(\Diff^+(\sphere^4))=G_{\mathbb{S}^5}^+(D^5,D^5)=\Upsilon_{\varphi}^+=\Gamma_{\varphi}^+\supset\{[\varphi]\},$$
$$ |G_{\mathbb{S}^5}(D^5,D^5)/\beta|=\Big\lfloor\frac{|\pi_0(\Diff^+(\sphere^4))|}{2}\Big\rfloor+1, \footnote{The supscript ``+" means that diffeomorphisms are orientation-preserving.}$$
which lead to
\begin{thm}\label{twopts focal}(\cite{Ge14})
\begin{itemize}
\item[1)] Each $\sphere^n$ $(n\neq 5)$ admits a unique equivalence class of isoparametric foliations with two points as focal submanifolds.
\item[2)] $\sphere^5$ admits a unique equivalence class of isoparametric foliations with two points as focal submanifolds if and only if $\pi_0(\Diff(\sphere^4))\simeq\mathbb{Z}_2$ (i.e., pseudo-isotopy implies isotopy on $\sphere^4$ as in other dimensions).
\end{itemize}
\end{thm}
We remark that so far no other information about $\pi_0(\Diff(\sphere^4))$ is known in literature.

For Step 1, \cite{Ge14} also presented many new examples. There exist many ``exotic"\footnote{Here by ``exotic" we mean non-equivalent as disk bundles but with diffeomorphic total spaces.} disk bundles $\widetilde{E}_{\pm}\cong E_{\pm}$. Hence there are many non-equivalent isoparametric foliations $(\widetilde{E}_{\varphi}, \widetilde{\fol}_{\varphi})\ncong (E_{\varphi},\fol_{\varphi})$:
\begin{itemize}
\item[(a)] \quad $\sphere^m\times D^k$ admits non-trivial disk bundle structures for (cf. \cite{DW00}) $$(m,k)=(7,4),(8,4),(9,4),(11,4),(11,5),(11,6).$$ Therefore each of $\sphere^{11},\sphere^{12},\sphere^{13},\sphere^{15},\sphere^{16},\sphere^{17}$ admits non-equivalent isoparametric foliations whose focal submanifolds are all standard spheres $(\sphere^m,\sphere^{k-1})$.
\item[(b)]\quad The tangent bundles of homotopy spheres $T\Sigma^n\cong T\sphere^n$ are ``exotic". Therefore, for instance, $\sphere^{14}$ admits $15$ (ignore orientation) non-equivalent isoparametric foliations whose focal submanifolds are $(\Sigma^7,\sphere^6)$ for $\Sigma^7\in\Theta_7\cong\mathbb{Z}_{28}$.
    \end{itemize}

\section{Application to existence of exotic smooth structures}
We start with the definitions of two inertia groups. Let $M^n$ be a closed oriented manifold. We will now recall the definition of the inertia subgroups $I_0(M)\subset \Theta_n$ and $I_1(M)\subset \Theta_{n+1}$ (cf. \cite{Le70}).

$I_0(M)$ consists of all $\Sigma\in \Theta_n$ such that $M\#\Sigma\cong M$. For
$\Sigma\in \Theta_n\backslash I_0(M)$, it is evident that $M\#\Sigma\ncong M$ is homeomorphic to $M$ and hence induces an exotic oriented smooth structure on $M$.   Moreover, different cosets in $\Theta_n/I_0(M)$ give distinct oriented smooth structures. Hence, there exist at least $|\Theta_n|/|I_0(M)|$ distinct oriented smooth structures on $M$.

To define $I_1(M)$, we first recall the disk theorem in \cite{Pa60}, i.e., any orientation-preserving diffeomorphism is isotopic to one that restricts to the identity on an embedded disk. Thus, for any $\phi\in\mathrm{Diff}^+(\sphere^{n})$, we can assume $\phi: \sphere^{n}=D^{n}_+\cup_{id}D^n_-\rightarrow D^{n}_+\cup_{id}D^n_-$ satisfies $\phi|_{D^{n}_+}=id$ up to isotopy. Then $I_1(M)$ consists of all $\Sigma_{\phi}\in\Theta_{n+1}$ such that the diffeomorphism of $M$ which differs from identity only on an $n$-disk in $M$, and there coincides with
$\phi$, is concordant to the identity. The coset space $\Theta_{n+1}/I_1(M)$ corresponds to a subset of $|\Theta_{n+1}|/|I_1(M)|$ elements in $\Gamma(M)$ (also in $\pi_0(\mathrm{Diff}^+(M))$).

In \cite{Le70}, Levine established an elegant relation between two inertia groups, i.e., $I_1(M)=I_0(M\times \sphere^1)$. Inspired by study of
isoparametric foliations on homotopy spheres, the following observation on inertia groups is acquired in \cite{Ge14}.
\begin{thm}(\cite{Ge14})
Let $M^{n-1}$ be a closed oriented embedded hypersurface in a closed oriented manifold $N^n$. Then $I_1(M^{n-1})\subseteq I_0(N^n)$. Consequently, $I_0(M^{n-1}\times\sphere^1)\subseteq I_0(N^n)$, i.e.,
$M^{n-1}\times\sphere^1$ has the smallest $I_0$ among all $n$-manifolds containing $M$.
\end{thm}
As an application of the theorem above, it follows from the fact the product of standard spheres has $I_0=0$ (cf. \cite{Sc71}) that
\begin{thm}(\cite{Ge14})
Let $M^{n-1}$ be a closed embedded hypersurface in $\sphere^n$. Then there exist at least $|\Theta_{n+k}|$ distinct oriented smooth structures on $M^{n-1}\times P^k\times\sphere^1$, where $P^k=\sphere^{k_1}\times\sphere^{k_2}\times\cdots \times \sphere^{k_l}$ is a product of standard spheres of total dimension $k=\sum_{i=1}^{l}k_i\geq 0$.
\end{thm}

\begin{ack}
We would like to thank the referees for their useful comments.
\end{ack}


\begin{thebibliography}{123}
\bibitem{AP2}
A. Akhmedov and B. Doug Park, \emph{Exotic smooth structures on small $4$-manifolds with odd signatures}, Invent. Math. \textbf{181} (2010), no. 3,
483--492.

\bibitem{ABT13}
M. M. Alexandrino, R. Briquet, and D. T\"oben, \emph{Progress in the theory of singular Riemannian foliations}, Differential Geom. Appl. \textbf{31}
(2013), no. 2, 248--267.

\bibitem{Ce08}
T. E. Cecil, \emph{Isoparametric and Dupin hypersurfaces}, SIGMA
\textbf{4} (2008), Paper 062, 28 pages.

\bibitem{CCJ07}
T. E. Cecil and Q. S. Chi and G. R. Jensen, \emph{Isoparametric
hypersurfaces with four principal curvatures}, Ann. Math.
\textbf{166} (2007), no.~1, 1--76.

\bibitem{Ce70}
 J. Cerf, \emph{La stratification naturelle des espaces de fonctions diff\'{e}rentiables r\'{e}elles et le th\'{e}or\`{e}me de la pseudo-isotopie}, Inst. Hautes \'{E}tudes Sci. Publ. Math. No. \textbf{39} (1970), 5--173.

\bibitem{Ch12}
Q. S. Chi, \emph{The isoparametric story}, National Taiwan University, June 25-July 6, 2012, \emph{http://www.math.wustl.edu/~chi/SummerCourse.pdf}

\bibitem{Ch13}
Q. S. Chi, \emph{Isoparametric hypersurfaces with four principal curvatures, III}, J. Diff. Geom. \textbf{94} (2013), 487--504.

\bibitem{DW00}
R. De Sapio and G. Walschap, \emph{Diffeomorphism of total spaces and equivalence of bundles}, Topology \textbf{39} (2000), 921--929.

\bibitem{DV12}
M. Dom\'{\i}nguez-V\'{a}zquez, \emph{Isoparametric foliations on complex projective spaces}, Trans. Amer. Math. Soc., DOI: http://dx.doi.org/10.1090/S0002-9947-2014-06415-5

\bibitem{DVG15}
M. Dom\'{\i}nguez-V\'{a}zquez and C. Gorodski, \emph{Polar foliations on quaternionic projective spaces}, preprint, 2015, arXiv:1507.02720.

\bibitem{Fi77}
R. Fintushel, \emph{Circle actions on simply connected $4$-manifolds}, Trans. Amer. Math.
Soc. \textbf{230} (1977), 147--171.

\bibitem{Fi78}
R. Fintushel, \emph{Classification of circle actions on $4$-manifolds}, Trans. Amer. Math. Soc.
\textbf{242} (1978), 377--390.

\bibitem{FQ90}
M. Freedman and F. Quinn, \emph{Topology of $4$-manifolds}, Princeton
Mathematical Series 39, Princeton University Press, Princeton, N.J.,
(1990).

\bibitem{GaRa13}
F. Galaz-Garcia and M. Radeschi, \emph{Singular Riemannian foliations and applications to positive and nonnegative curvature}, J. Topology \textbf{8} (2015), 603--620.

\bibitem{Ge14}
J. Q. Ge, \emph{Isoparametric foliations, diffeomorphism groups and exotic smooth structures}, preprint, 2014, arXiv:1404.6194.

\bibitem{GR13}
J. Q. Ge and M. Radeschi, \emph{Differentiable classification of $4$-manifolds with singular Riemannian foliations}, Math. Ann. \textbf{363} (2015) 525--548.

\bibitem{GT12}
J. Q. Ge and Z. Z. Tang, \emph{Chern conjecture and isoparametric hypersurfaces}, in ¡°Differential
Geometry - under the influence of S.S.Chern¡±. Edited by Y. B. Shen, Z.M. Shen and S. T. Yau,
Higher Education Press and International Press, Beijing-Boston, 2012.

\bibitem{GT13}
J. Q. Ge and Z. Z. Tang, \emph{Isoparametric functions and exotic
spheres}, J. Reine Angew. Math. \textbf{683} (2013), 161-180.

\bibitem{GTY14} J. Q. Ge, Z. Z. Tang and W. J. Yan, \emph{A filtration for isoparametric
hypersurfaces in Riemannian manifolds}, J. Math. Soc. Japan. \textbf{67} (2015), 1179--1212.

\bibitem{GW09}
D. Gromoll and G. Walschap, \emph{Metric foliations and curvature}, Progress in Mathematics, \textbf{268}, Birkh\"{a}user Verlag, Basel, 2009.

\bibitem{GH87}
K. Grove and S. Halperin, \emph{Dupin hypersurfaces, group actions and the double mapping cylinder}, J. Diff. Geom. \textbf{26} (1987), 429--459.

\bibitem{GW14}
K. Grove and B. Wilking, \emph{A knot characterization and 1-connected nonnegatively curved 4-manifolds with circle symmetry}, Geom. Topol. \textbf{18} (2014), 3091--3110.

\bibitem{GZ00}
K. Grove and W. Ziller, \emph{Curvature and Symmetry of Milnor
Spheres}, Ann. Math. \textbf{152} (2000), 331--367.

\bibitem{Ho10}
C. A. Hoelscher, \emph{Classification of cohomogeneity one manifolds in low dimensions}, Pac. J. Math. \textbf{246} (2010), 129--185.

\bibitem{JW08}
M. Joachim and D. J. Wraith, \emph{Exotic spheres and curvature}, Bull. Amer. Math. Soc. (N.S.) \textbf{45} (2008), no. 4, 595--616.

\bibitem{KM63} M. A. Kervaire and J. Milnor, \emph{Group of homotopy spheres: I}, Ann. Math. \textbf{77} (1963), 504--537.

\bibitem{Le70}
J. Levine, \emph{Inertia groups of manifolds and diffeomorphisms of spheres}, Amer. J. Math. \textbf{92} (1970), 243--258.

\bibitem{LW13}
A. Lytchak, B. Wilking, \emph{Riemannian foliations of spheres}, preprint, 2013,  arXiv:1309.7884.

\bibitem{Mil56}
J. Milnor, \emph{On manifolds homeomorphic to the $7$-sphere}, Ann.
Math. (2) \textbf{64} (1956), 399--405.

\bibitem{Mi13}
R. Miyaoka, \emph{Isoparametric hypersurfaces with $(g,m)=(6,2)$}, Ann. Math. \textbf{177} (2013), 53--110.

\bibitem{Mol88}
P. Molino, \emph{Riemannian foliations}, Translated from the French by Grant Cairns. With appendices by G. Cairns, Y. Carri\`{e}re, E. Ghys, E. Salem,
V. Sergiescu. Progress in Math. \textbf{73}. Birkh\"{a}user Boston, Inc., Boston, MA, 1988.

\bibitem{OR70}
P. Orlik and F. Raymond, \emph{Actions of the torus on $4$-manifolds I}, Trans. Amer. Math.
Soc. \textbf{152} (1970), 531--559.

\bibitem{OR74}
P. Orlik and F. Raymond, \emph{Actions of the torus on $4$-manifolds. II}, Topology \textbf{13}
(1974), 89--112.

\bibitem{Pa60}
R. S. Palais, \emph{Extending diffeomorphisms}, Proc. Amer. Math. Soc. \textbf{11} (1960), 274--277.

\bibitem{Pa86}
J. Parker, \emph{$4$-dimensional G-manifolds with $3$-dimensional orbits}, Pac. J. Math. \textbf{125} (1986), 187--204.

\bibitem{QT13}
C. Qian and Z. Z. Tang, \emph{Isoparametric functions on exotic spheres}, Advances in Math. \textbf{272} (2015), 611--629.


\bibitem{QT14}
C. Qian and Z. Z. Tang, \emph{Recent progress in isoparametric functions and isoparametric hypersurfaces}, Real and Complex Submanifolds, Springer Proceedings in Mathematics \& Statistics \textbf{106} (2014), 65--76.

\bibitem{Ra14}
M. Radeschi, \emph{Clifford algebras and new singular Riemannian foliations in spheres},  Geom. Funct. Anal.  \textbf{24} (2014), 1660--1682.

\bibitem{Sc71}
R. Schultz, \emph{On the inertia group of a product of spheres}, Trans. Amer. Math. Soc. \textbf{156} (1971), 137--153.

\bibitem{TXY14}
Z. Z. Tang, Y. Q. Xie and W. J. Yan,  \emph{Isoparametric foliation and Yau conjecture on the first eigenvalue,II}, J. Funct. Anal. \textbf{266}, (2014), 6174--6199.

\bibitem{TY13}
Z. Z. Tang and W. J. Yan,  \emph{Isoparametric foliation and Yau conjecture on the first eigenvalue}, J. Diff. Geom. \textbf{94} (2013), 521--540.

\bibitem{TY15}
Z. Z. Tang and W. J. Yan, \emph{Isoparametric foliation and a problem of Besse on generalizations of Einstein condition}, Advances in Math. \textbf{285} (2015), 1970--2000.

\bibitem{Th91}
G. Thorbergsson, \emph{Isoparametric foliations and their buildings}, Ann. Math. (2)  \textbf{133} (1991),  no. 2, 429--446.

\bibitem{Th00}
G. Thorbergsson, \emph{A survey on isoparametric hypersurfaces and
their generalizations}, In Handbook of differential geometry, Vol.
I, North - Holland, Amsterdam, (2000), pages 963 - 995.

\bibitem{Th10}
G. Thorbergsson, \emph{Singular Riemannian foliations and isoparametric submanifolds}, Milan J. Math. \textbf{78} (2010), 355--370.

\bibitem{Wa87}
Q. M. Wang, \emph{Isoparametric Functions on Riemannian Manifolds. I}, Math. Ann. \textbf{277} (1987), 639--646.

\end{thebibliography}
\end{document}